\documentclass[11pt]{amsart}
\usepackage{amsfonts}
\usepackage{amsmath}
\usepackage{amssymb}
\usepackage[mathscr]{eucal}
\usepackage{url}

\newcommand{\ph}[2]{{\left({#1}\right)}_{#2}}

\newcommand{\gf}[1]{\Gamma{\left({#1}\right)}}

\renewcommand*{\bar}{\overline}

\newcommand{\bin}[2]{\left({\genfrac{}{}{0pt}{}{#1}{#2}}\right)}

\theoremstyle{plain}
\newtheorem{theorem}{Theorem}[section]

\newtheorem{prop}[theorem]{Proposition}
\newtheorem{cor}[theorem]{Corollary}
\theoremstyle{definition}

\numberwithin{equation}{section}

\begin{document}

\title[$_3F_2$ hypergeometric series and periods of elliptic curves]{$_3F_2$ hypergeometric series and periods of elliptic curves}
\author{Dermot M\lowercase{c}Carthy}

\address{School of Mathematical Sciences, University College Dublin, Belfield, Dublin 4, Ireland}

\email{dermot.mc-carthy@ucdconnect.ie}


\date{September 3, 2008}

\begin{abstract}
We express the real period of a family of elliptic curves in terms of classical hypergeometric series. This expression is analogous to a result of Ono which relates the trace of Frobenius of the same family of elliptic curves to a Gaussian hypergeometric series. This analogy provides further evidence of the interplay between classical and Gaussian hypergeometric series. 
\end{abstract}

\maketitle

\section{Introduction}

In \cite{G}, Greene introduced the notion of general hypergeometric series over finite fields or \emph{Gaussian hypergeometric series}, which are analogous to classical hypergeometric series. The motivation for his work was to develop the area of character sums and their evaluations through parallels with the theory of hypergeometric functions. The basis for this parallel was the analogy between Gauss sums and the gamma function as discussed in \cite{Ev, IR, Ko3, Y}.

Since then, the interplay between ordinary hypergeometric series and Gaussian hypergeometric series has played an important role in character sum evaluations \cite{GS}, supercongruences \cite{M}, finite field versions of the Lagrange inversion formula \cite{G2} and the representation theory of SL(2, $\mathbb{R}$) \cite{G3}. Recently, the author in \cite{Ro} has further developed this interplay by providing an expression for the real period of an elliptic curve in Legendre normal form in terms of an ordinary hypergeometric series. This formula is analogous to an expression for the trace of Frobenius of the curve in terms of a Gaussian hypergeometric series. He then displays a striking analogy between binomial coefficients involving rational numbers and those involving multiplicative characters. This paper examines this analogy further using a different family of elliptic curves and is organized as follows. In Section 2 we outline this analogy and state our results. Section 3 recalls some properties of ordinary hypergeometric series, elliptic curves and the arithmetic-geometric mean. In Section 4 we prove our results.

\section{Statement of Results}
We recall that the ordinary hypergeometric series $_pF_q$ is defined by
\begin{equation*}
{_pF_q} \left( \begin{array}{ccccc} a_1, & a_2, & a_3, & \dotsc, & a_p \vspace{.05in}\\
\phantom{a_1} & b_1, & b_2, & \dotsc, & b_q \end{array}
\Big| \; z \right)
:=\sum^{\infty}_{n=0}
\frac{\ph{a_1}{n} \ph{a_2}{n} \ph{a_3}{n} \dotsm \ph{a_p}{n}}
{\ph{b_1}{n} \ph{b_2}{n} \dotsm \ph{b_q}{n}}
\; \frac{z^n}{{n!}}
\end{equation*}
where $a_i$, $b_i$ and $z$ are complex numbers, with none of the $b_i$ being negative integers or zero, $p$ and $q$ are positive integers, $\ph{a}{0}:=1$ and $\ph{a}{n} := a(a+1)(a+2)\dotsm(a+n-1)$ for positive integers $n$.

Let $\mathbb{F}_{p}$ denote the finite field with $p$ elements. We extend the domain of all characters $\chi$ of $\mathbb{F}^{*}_{p}$ to $\mathbb{F}_{p}$, by defining $\chi(0):=0$. We now introduce two definitions from \cite{G}. The first definition is the finite field analogue of the binomial coefficient. For characters $A$ and $B$ of $\mathbb{F}_{p}$, define $\bin{A}{B}$ by
\begin{equation}\label{FF_Binomial}
\binom{A}{B} := \frac{B(-1)}{p} J(A, \bar{B})
\end{equation}
where $J(\chi, \lambda)$  denotes the Jacobi sum for $\chi$ and $\lambda$ characters of $\mathbb{F}_{p}$. The second definition is the finite field analogue of ordinary hypergeometric series. For characters $A_0,A_1,\dotsc, A_n$ and $B_1, \dotsc, B_n$ of $\mathbb{F}_{p}$ and 
$x \in \mathbb{F}_{p}$, define the \textit{Gaussian hypergeometric series} by
\begin{equation*}
{_{n+1}F_n} {\left( \begin{array}{cccc} A_0, & A_1, & \dotsc, & A_n \\
\phantom{A_0} & B_1, & \dotsc, & B_n \end{array}
\Big| \; x \right)}_{p}
:= \frac{p}{p-1} \sum_{\chi} \binom{A_0 \chi}{\chi} \binom{A_1 \chi}{B_1 \chi}
\dotsm \binom{A_n \chi}{B_n \chi} \chi(x)
\end{equation*}
where the summation is over all characters $\chi$ on $\mathbb{F}_{p}$. In the case where $A_i = \phi_p$, the quadratic character mod $p$, for all $i$ and $B_j= \epsilon_p$, the trivial character mod $p$, for all $j$ we denote this by ${_{n+1}F_n}(x)_{p}$ for brevity. 

We now briefly recall some facts about elliptic curves. For further details see \cite{Kn}, \cite{Ko2} and \cite{S}. Recall that every elliptic curve $E / \mathbb{C}$ can be written in the form
\begin{equation}\label{Egg}
y^2 = 4x^3 - g_2x - g_3 \: ,
\end{equation}
with $g_2$, $g_3 \in \mathbb{C}.$ We can associate a period lattice $\Lambda$ to $E$ via the biholomorphic mapping $\varphi$ :  $\mathbb{C} / \Lambda \rightarrow E(\mathbb{C}) $ given by
\begin{equation*}
\varphi(z) = 
\left\{ \begin{array}{ll} (\wp(z), \wp'(z), 1) & \quad \textup{for } z \notin \Lambda, \vspace{.05in} \\
(0,1,0) & \quad \textup{for } z \in \Lambda, \end{array} \right.
\end{equation*}
where $\wp$ is the Weierstrass $\wp$-function. If $g_2$, $g_3 \in \mathbb{R}$ then $\Lambda$ can be chosen to be of the form $\Lambda = \Omega(E) \mathbb{Z} + \Omega'(E) \mathbb{Z}$ where $\Omega(E) \in \mathbb{R}$ and $\Omega'(E) \in \mathbb{C}$. We call $\Omega(E)$ the \emph{real period} of $E$. Furthermore, if the right-hand side of (\ref{Egg}) has three real roots then $\Omega'(E)$ will be strictly imaginary.

Consider an elliptic curve $E / \mathbb{Q}$ in Weierstrass form
\begin{equation*}
E: y^2 + a_1 xy + a_3y = x^3 + a_2 x^2 +a_4 x + a_6
\end{equation*}
where $a_i \in \mathbb{Q}$. Defining the quantities
\begin{equation*}
b_2 := {a_1}^2 +4 a_2, b_4 := 2 a_4 + a_1a_3, b_6 := {a_3}^2 + 4 a_6,
\end{equation*}
and
\begin{equation*}
b_8 := {a_1}^2 a_6 + 4 a_2 a_6 - a_1 a_3 a_4 + a_2 {a_3}^2 - {a_4}^2,
\end{equation*}
the discriminant of $E$, $\Delta(E)$, is given by
\begin{equation*}
\Delta(E) = -{b_2}^2 b_8 - 8{b_4}^3 - 27{b_6}^2 + 9 b_2 b_4 b_6.
\end{equation*}
Let $\tilde{E}$ denote the reduction of $E$ mod $p$. Recall that if $p \nmid \Delta(E)$ then $E$ has good reduction ($\tilde{E} / \mathbb{Q}$ is an elliptic curve) and we say $p$ is a prime of good reduction. We define the integer $a_p(E)$ by
\begin{equation*}
a_p(E) := 1 + p - N_p \: ,
\end{equation*}
where $N_p$ is the number of rational points on $\tilde{E}$ over $\mathbb{F}_p$ (including the point at infinity). When $p$ is a prime of good reduction, we refer to $a_p(E)$ as the \emph{trace of Frobenius} as it can be interpreted as the trace of the Frobenius endomorphism on $E$. Furthermore, if $E$ is given by $y^2 = f(x)$ then
\begin{equation}\label{TraceFormula}
a_p(E) = - \sum_{x \in \mathbb{F}_p} \phi_p(f(x)) \: .
\end{equation}

Consider the family of elliptic curves $E_{\lambda} / \mathbb{Q}$ defined by
\begin{equation*}
E_{\lambda}: y^2=(x-1)(x^2+\lambda), \qquad \lambda \in \mathbb{Q} \setminus \lbrace 0, -1 \rbrace \: .
\end{equation*}
\noindent Ono \cite[Thm. 5]{O}, (see also \cite[Chapter 11]{O2}), proved that if $\lambda \in \mathbb{Q} \setminus \lbrace 0,-1\rbrace$ and $p$ is an odd prime for which ord$_p(\lambda(\lambda+1)) = 0$ then
\begin{equation}\label{ThmOno}
_3F_2 \left( \frac{1+\lambda}{\lambda} \right)_p
= \frac{\phi_p(-\lambda)\left(a_p(E_\lambda)^2 - p\right)}{p^2} \: .
 \end{equation}
Note that a change of variables in Theorem 5 of \cite{O} is required to arrive at (\ref{ThmOno}) (see also \cite[Thm. 4.4]{FOP}). The condition ord$_p(\lambda(\lambda+1)) = 0$ ensures that $p$ is a prime of good reduction and so $a_p(E_\lambda)$ is the trace of Frobenius. Using the following property of Gaussian hypergeometric series (see \cite[Thm. 4.2]{G})
\begin{equation*}
_3F_2\left(\frac{1}{t}\right)_p = \phi_p(-t) _3F_2(t)_p \: ,
\end{equation*}
we transform (\ref{ThmOno}) to get
\begin{equation*}
_3F_2 \left( \frac{\lambda}{1+\lambda} \right)_p
= \frac{\phi_p(1+\lambda)\left(a_p(E_\lambda)^2 - p \right)}{p^2} \: .
 \end{equation*}
 
We would like to prove an analogous formula which replaces the Gaussian hypergeometric series with the classical hypergeometric series and trace of Frobenius with the real period $\Omega(E_\lambda)$. In \cite{Ro} this analogy is based on replacing characters of order $n$ in the Gaussian hypergeometric series with $1/n$ as the arguments in the classical hypergeometric series, $\phi_p(s)$ with $\sqrt{s}$ and $\phi_p(-1) p$ with $\pi$. This would suggest
\begin{equation}\label{Guess}
{_3F_2} \left( \begin{array}{ccc} \frac{1}{2}, & \frac{1}{2}, & \frac{1}{2} \\
\phantom{\frac{1}{2},} & 1, & 1 \end{array}
\Big| \; \frac{\lambda}{1+\lambda} \right) =
\frac{\sqrt{1+\lambda} \: \:  \Omega(E_\lambda)^2}{\pi^2} - i \frac{\sqrt{1+\lambda}}{\pi}
\end{equation}
as an appropriate analogy. The main result of this paper extends the analogy in \cite{Ro} by taking the real part of the right-hand side of (\ref{Guess}). This extension is consistent with the results in \cite{Ro}. 
\begin{theorem}\label{TheoremPeriod}
Let $E_\lambda$ be the elliptic curve defined by
\begin{equation*}
E_{\lambda}: y^2=(x-1)(x^2+\lambda), \qquad \lambda \in \mathbb{R} \setminus \lbrace 0, -1 \rbrace \: .
\end{equation*}
Then for $\lambda>0$,
\begin{equation*}
{_3F_2} \left( \begin{array}{ccc} \frac{1}{2}, & \frac{1}{2}, & \frac{1}{2} \\
\phantom{\frac{1}{2},} & 1, & 1 \end{array}
\Big| \; \frac{\lambda}{1+\lambda} \right) =
\frac{\sqrt{1+\lambda} \: \:  \Omega(E_\lambda)^2}{\pi^2} \: ,
\end{equation*}
where $\Omega(E_\lambda)$ is the real period of $E_\lambda$.
\end{theorem}
\noindent Note that the analogy in \cite{Ro} also contained a factor of ${-1}$ which the author explains is inherent in the definition of Gaussian hypergeometric series. This may be better explained by the $-1$ preceding the character sum expression for $a_p(E_\lambda)$ in (\ref{TraceFormula}), which disappears upon squaring. (The real period can be expressed as an elliptic integral which is somewhat analogous to (\ref{TraceFormula}) but without the minus sign). Therefore we suggest a further refinement which would see $a_p(E_\lambda)$ being replaced with $-\Omega(E_\lambda)$.

We now specialize the curve by choosing $\lambda =$ 1/3. We then use a known transformation of the hypergeometric series in terms of the gamma function to simplify the expression as a binomial coefficient. We extend the interpretation of the binomial coefficient to include rational arguments via
\begin{equation*}
\bin{n}{k} = \frac{\gf{n+1}}{\gf{k+1}\gf{n-k+1}} \: .
\end{equation*}
Our result is as follows.
 
\begin{cor}\label{CorBinomial}
Let $E_{\frac{1}{3}}$ be the elliptic curve defined by
\begin{equation*}
E_{\frac{1}{3}}: y^2=(x-1)(x^2+\tfrac{1}{3}) \: .
\end{equation*}
Then
\begin{equation*}
\frac{2 \sqrt{2}}{3\pi} \cdot \Omega(E_\frac{1}{3}) = \bin{\frac{1}{3}}{\frac{1}{2}}
\end{equation*}
and 
\begin{equation*}
\sqrt{2} \cdot \Omega(E_\frac{1}{3}) = \frac{\gf{\frac{1}{3}} \gf{\frac{1}{2}}}{\gf{\frac{5}{6}}} \: ,
\end{equation*}
where $\Omega(E_\frac{1}{3})$ is the real period of $E_{\frac{1}{3}}$.
\end{cor}

\noindent We now find an analogous result in terms of the trace of Frobenius and the binomial coefficient of characters, as defined in (\ref{FF_Binomial}), which we can also express in terms of Gauss sums.

\begin{theorem}\label{TheoremBinomial}
Let $E_{\frac{1}{3}}$ be the elliptic curve defined by
\begin{equation*}
E_{\frac{1}{3}}: y^2=(x-1)(x^2+\tfrac{1}{3}) \: .
\end{equation*}
Then for $p$ a prime with $p>3$,
\begin{equation}\label{ThmBinPart1}
- \frac{\phi_p(-2)}{p} \cdot a_p(E_\frac{1}{3}) = 2 \: \textup{Re} \bin{\chi_3}{\phi_p}
\end{equation}
and
\begin{equation}\label{ThmBinPart2}
- \phi_p(2) \cdot  a_p(E_\frac{1}{3}) = 2 \: \textup{Re} \left[\frac{G(\chi_3) G(\phi_p)}{G(\chi_3 \: \phi_p)}\right] \: ,
\end{equation}
where $a_p(E_\frac{1}{3})$ is the trace of Frobenius of $E_{\frac{1}{3}}$, $\chi_3$ is a character of order three of $\mathbb{F}_p$ and $G(\chi)$ is a Gauss sum.
\end{theorem}

\noindent Again the analogy is achieved by replacing $\Omega(E_\frac{1}{3})$ with $-a_p(E_\frac{1}{3})$, $\sqrt{s}$ with $\phi_p(s)$, $\pi$ with $\phi_p(-1) p$, rational numbers with characters and taking the real part of the terms involving characters. We now also replace the gamma function with the Gauss sum in the second result, which is what we would expect. The analogy holds up to a factor of 2 (which also appears in \cite{Ro}). This might be explained as follows. It should be possible to express the trace of Frobenius in terms of a $_2F_1$ Gaussian hypergeometric series which in turn could be expressed as a Jacobi sum plus its conjugate, using Theorem 4.16 in \cite{G}, which would evaluate as two times its real part. However, we do not investigate this here. 

\section{Preliminaries}
We first recall some properties of ordinary hypergeometric series. In particular we recall that a $_2F_1$ has the following integral representation \cite[page 115 (7)]{E}:
\begin{equation}\label{IntegralRep}
{_2F_1} \left( \begin{array}{cc} a, & b \\
\phantom{a,} & c \end{array}
\Big| \; z \right)  = \frac{2 \: \gf{c}}{\gf{b} \gf{c-b}} \int_{0}^{\frac{\pi}{2}} \frac{(\sin{t})^{2b-1} (\cos{t})^{2c-2b-1}}{(1-z \sin^2{t})^a} dt 
\end{equation}
where Re $c >$ Re $b > 0$. We also note three transformation properties which we will use in Section 4.
From \cite[page 111 (10)]{E} we have that
\begin{equation}\label{TransformQuadratic}
{_2F_1} \left( \begin{array}{cc} a, & b \\
\phantom{a,} & a+b+\tfrac{1}{2} \end{array}
\Big| \; z \right)  = 
{_2F_1} \left( \begin{array}{cc} 2a, & 2b \\
\phantom{2a,} & a+b+\tfrac{1}{2} \end{array}
\Big| \; \tfrac{1}{2} - \tfrac{1}{2}(1-z)^{\frac{1}{2}} \right) \: ,
\end{equation}
from \cite[Entry 33(iii)]{B} that
\begin{equation}\label{TransformSquare}
{_3F_2} \left( \begin{array}{ccc} \frac{1}{2}, & \frac{1}{2}, & \frac{1}{2} \\
\phantom{\frac{1}{2},} & 1, & 1 \end{array}
\Big| \; z \right)  =
{\left[{_2F_1} \left( \begin{array}{cc} \frac{1}{4}, & \frac{1}{4} \\
\phantom{\frac{1}{4},} & 1 \end{array}
\Big| \; z \right) \right]}^2,
\end{equation}
and from \cite[page 105 (3)]{E} that
\begin{equation}\label{Transform3}
{_2F_1} \left( \begin{array}{cc} a, & b \\
\phantom{a,} & c \end{array}
\Big| \; z \right)  = (1-z)^{-a}
{_2F_1} \left( \begin{array}{cc} a, & c-b \\
\phantom{a,} & c \end{array}
\Big| \; \frac{z}{z-1} \right).
\end{equation}
These transformations are valid for all values of $z$ for which the series involved converge.

We now recall the definition of the Arithmetic-Geometric Mean. Given two positive real numbers $\alpha$ and $\beta$, the Arithmetic-Geometric Mean of $\alpha$ and $\beta$, denoted AGM($\alpha$,$\beta$), is defined as the common limit of the two sequences $\alpha_n$ and $\beta_n$ where $\alpha_0:=\alpha$, $\beta_0:=\beta$, $\alpha_{n+1}:=(\alpha_n+\beta_n)/2$ and $\beta_{n+1}:=\sqrt{\alpha_n \beta_n}$.
The AGM can be expressed as an integral (\cite[page 390]{Co}) which can be transformed into a hypergeometric series as follows.
\begin{align}\label{AGM_to_Hyp}
\frac{\pi}{AGM(\alpha, \beta)}
&= 2 \; \int_{0}^{\frac{\pi}{2}} \frac{dt} {\sqrt{\alpha^2 \cos^2{t} +  \beta^2 \sin^2{t}}}\\
&\notag= 2 \alpha^{-1}  \int_{0}^{\frac{\pi}{2}} \left[{ \cos^2{t} +  {\left(\tfrac{\beta}{\alpha}\right)}^2 \sin^2{t} }\right]^{-\frac{1}{2}} dt\\
&\notag= 2 \alpha^{-1}  \int_{0}^{\frac{\pi}{2}}  \left[{{1-  \left(1-{\left(\tfrac{\beta}{\alpha}\right)}^2\right) \sin^2{t}}}\right] ^{-\frac{1}{2}} dt \\
&\notag= \alpha^{-1} \: \pi \; {_2F_1} \left( \begin{array}{cc} \frac{1}{2}, & \frac{1}{2} \\
\phantom{\frac{1}{2},} & 1 \end{array}
\Big| \; 1- \left(\tfrac{\beta}{\alpha}\right)^2 \right) \: .
\end{align}
The last step in (\ref{AGM_to_Hyp}) follows from (\ref{IntegralRep}) with $a=b=\frac{1}{2}$, $c=1$ and $z=1-{\left(\tfrac{\beta}{\alpha}\right)}^2$.

Next we introduce the notion of a \emph{quadratic twist} of an elliptic curve. Let $E /\mathbb{Q}$ be an elliptic curve defined by
\begin{equation*}
E: y^2 = x^3 +ax^2 + bx +c,
\end{equation*}
with $a,b,c \in \mathbb{Q}$. If $t$ is a square-free integer, then the $t$-quadratic twist of $E$, which we denote $E_t$, is defined by
\begin{equation*}
E_t : y^2 = x^3 +atx^2 + bt^2x + ct^3.
\end{equation*}
If $p$ is a prime of good reduction for both $E$ and $E_t$ and gcd($p$, 6)=1, then
\begin{equation}\label{Twist}
 a_p(E) =  \phi_p(t) \: a_p(E_t).
 \end{equation}

We now mention a result which we will use in the proof of Theorem \ref{TheoremBinomial}. As it follows from well-known properties of Jacobsthal sums (see Sections 6.1 and 6.2 in \cite{BEW}) we omit the proof.
\begin{prop}\label{PropJacobsthal}
\begin{equation*}
\sum_{x \in \mathbb{F}_p} \phi_p(x^3+1)
= \left\{ \begin{array}{ll} 2a & \quad \textup{if} \quad p\equiv1 \pmod 3, \textup{ where } p=a^2+3b^2 \textup{ and } a\equiv-1 \pmod 3, \vspace{.05in} \\
\phantom{2}0 & \quad \textup{if} \quad p\equiv2 \pmod 3 \; .\end{array} \right.
\end{equation*}
\end{prop}

\section{Proofs of Theorem \ref{TheoremPeriod}, Corollary \ref{CorBinomial} and Theorem \ref{TheoremBinomial}  }
\begin{proof}[Proof of Theorem \ref{TheoremPeriod}]
Making the change of variable $y \mapsto \frac{y}{2}$ in $E_\lambda$ yields 
\begin{equation*}
E'_{\lambda}: y^2=4(x-1)(x^2+\lambda)\: .
\end{equation*} 
For $\lambda>0$, we note that $4(x-1)(x^2+\lambda)=0$ has one real root. The real period $\Omega(E_\lambda)$ is then given by \cite[page 391]{Co}
\begin{equation*}
\Omega(E_\lambda) = \frac{2\pi}{AGM(2\sqrt{b},\sqrt{2b+a})} \: ,
\end{equation*}
where
$a=2$ and $b=\sqrt{1+\lambda}$. Using (\ref{AGM_to_Hyp}) we get
\begin{equation*}
\Omega(E_\lambda) = {\left(\sqrt{1+\lambda}\right)}^{-\frac{1}{2}} \; \pi \;
{_2F_1} \left( \begin{array}{cc} \frac{1}{2}, & \frac{1}{2} \\
\phantom{\frac{1}{2},} & 1 \end{array}
\Big| \; \tfrac{1}{2} \left(1- \tfrac{1}{\sqrt{1+\lambda}}\right) \right) \: .
\end{equation*}
Now applying equation (\ref{TransformQuadratic}) with $a=b=\frac{1}{4}$ and $z=\frac{\lambda}{1+\lambda}$ we get
\begin{equation}\label{Period2F1}
\Omega(E_\lambda) = {\left(\sqrt{1+\lambda}\right)}^{-\frac{1}{2}} \; \pi \;
{_2F_1} \left( \begin{array}{cc} \frac{1}{4}, & \frac{1}{4} \\
\phantom{\frac{1}{4},} & 1 \end{array}
\Big| \; \frac{\lambda}{1+\lambda} \right) \: .
\end{equation}
We note that the condition $\lambda>0$ implies that these hypergeometric series converge. Squaring both sides and applying (\ref{TransformSquare}) yields the result.
\end{proof}

\begin{proof}[Proof of Corollary \ref{CorBinomial}]
Transforming the hypergeometric series on the right-hand side of (\ref{Period2F1}) using (\ref{Transform3}), with $a=b=\frac{1}{4}$ and $z=\frac{\lambda}{1+\lambda}$, we see that for $0<\lambda<1$,
\begin{equation*}
\Omega(E_\lambda) = \pi \;
{_2F_1} \left( \begin{array}{cc} \frac{1}{4}, & \frac{3}{4} \\
\phantom{\frac{1}{4},} & 1 \end{array}
\Big| \; -\lambda \right) \:.
\end{equation*}
Now letting $\lambda=\frac{1}{3}$ and noting that
\begin{equation*}
{_2F_1} \left( \begin{array}{cc} \frac{1}{4}, & \frac{3}{4} \\
\phantom{\frac{1}{4},} & 1 \end{array}
\Big| \; -\tfrac{1}{3} \right)
= \frac{3}{2\sqrt{2}} \cdot \frac{\gf{\frac{4}{3}}}{\gf{\frac{3}{2}} \gf{\frac{5}{6}}}\: ,
\end{equation*}
(see \cite[page 104 (53)]{E} with $a=-\frac{1}{4}$), the first result follows.
The second result then follows upon recalling that $\gf{1+x} = x \: \gf{x}$ and ${\gf{\frac{1}{2}}}^2 = \pi$.
\end{proof}

\begin{proof}[Proof of Theorem \ref{TheoremBinomial}]
We first prove (\ref{ThmBinPart1}). By definition (\ref{FF_Binomial}) it suffices to prove
\begin{equation}\label{NewTheorem}
-\tfrac{1}{2} \cdot {\phi_p(2)} \cdot a_p(E_\frac{1}{3}) = \textup{Re} [J(\chi_3, \phi_p)] \: .
\end{equation}
We now evaluate $a_p(E_\frac{1}{3}$). A similar calculation appears in \cite{O} although we present our result slightly differently. Making the change of variables $(x,y) \mapsto (\frac{x}{9}+\frac{1}{3},\frac{y}{27})$ in $E_\frac{1}{3}$ yields
\begin{equation*}
{E'}_\frac{1}{3}: y^2 = x^3 - 6^3 \: .
\end{equation*}
This is the $-6$-quadratic twist of $y^2=x^3+1$. Therefore, applying (\ref{TraceFormula}) and (\ref{Twist}), noting that for $p>3$, gcd($p$, 6)=1 and $p$ is a prime of good reduction for both $y^2 = x^3 + 1$ and ${E'}_\frac{1}{3}$, we get
\begin{equation*}
 a_p(E_\frac{1}{3}) =  a_p({E'}_\frac{1}{3}) = - \phi_p(-6) \sum_{x \in \mathbb{F}_p} \phi_p(x^3+1) \: .
\end{equation*}
Using Proposition \ref{PropJacobsthal} and the fact that $\phi_p(-3) = 1$ if and only if $p \equiv 1\pmod3$ we get
\begin{equation*}
-\tfrac{1}{2} \cdot \phi_p(2) \cdot a_p(E_\frac{1}{3}) 
= \left\{ \begin{array}{ll} a & \quad \textup{if} \quad p\equiv1 \pmod 3, \textup{ where } p=a^2+3b^2 \textup{ and } a\equiv-1 \pmod 3, \vspace{.05in} \\
0 & \quad \textup{if} \quad p\equiv2 \pmod 3 \; .\end{array} \right.
\end{equation*}
Next we examine the right-hand side of (\ref{NewTheorem}). We first note that
\begin{equation*}
J(\chi_3, \phi_p) = \sum_{x \in \mathbb{F}_p} \chi_3(x) \phi_p(1-x) \: .
\end{equation*}
If $p \equiv 2 \pmod 3$ then $\chi_3(x) = 1$ for all $x \in \mathbb{F}_p$. Therefore,
\begin{equation*}
J(\chi_3, \phi_p) = \sum_{x \in \mathbb{F}_p} \phi_p(1-x) = 0 \;.
\end{equation*}
If $p \equiv 1 \pmod 3$, then
\begin{align*}
\textup{Re} [J(\chi_3, \phi_p)]
&= \textup{Re} \left[\sum_{x \in \mathbb{F}_p} \chi_3(x) \phi_p(1-x)\right]\\
\notag &= \sum_{x \in {\mathbb{F}_p^*}^3} \phi_p(1-x)
+ \textup{Re} \left[\sum_{x \in \mathbb{F}_p^* \setminus {\mathbb{F}_p^*}^3} \chi_3(x) \phi_p(1-x)\right]\\
\notag &= \frac{1}{3} \sum_{x \in \mathbb{F}_p^*} \phi_p(1-x^3)
+ \sum_{x \in \mathbb{F}_p^* \setminus {\mathbb{F}_p^*}^3} \textup{Re}\left[\chi_3(x)\right] \phi_p(1-x)\\
\notag &= \frac{1}{3} \sum_{x \in \mathbb{F}_p^*} \phi_p(1-x^3)
- \frac{1}{2} \sum_{x \in \mathbb{F}_p^* \setminus {\mathbb{F}_p^*}^3} \phi_p(1-x)\\
\notag &= \frac{1}{3} \sum_{x \in \mathbb{F}_p^*} \phi_p(1-x^3)
- \frac{1}{2} \left[-1- \frac{1}{3} \sum_{x \in \mathbb{F}_p^*} \phi_p(1-x^3) \right]\\
\notag &= \frac{1}{2} \sum_{x \in \mathbb{F}_p}\phi_p(1-x^3)\\
\notag &= \frac{1}{2} \sum_{x \in \mathbb{F}_p}\phi_p(x^3+1)\\
\notag &= a \qquad (\textup{by Proposition (\ref{PropJacobsthal})})
\end{align*}
where $p=a^2+3b^2$ and $a\equiv-1 \pmod 3$, which completes the proof of (\ref{ThmBinPart1}). Expanding the right-hand side of (\ref{ThmBinPart1}) using (\ref{FF_Binomial}), and then using the fact that
\begin{equation*}
J(\chi, \psi) = \frac{G(\chi) G(\psi)}{G(\chi \: \psi)},
\end{equation*}
where $J(\chi, \psi)$ and $G(\chi)$ are Jacobi and Gauss sums respectively, (\ref{ThmBinPart2}) follows.
\end{proof}

\section{Remark}
It is worth noting that if $E / \mathbb{R}$ is an elliptic curve and $E(\mathbb{R})$ has a single connected component then $E$ is isomorphic to either $E_\lambda$, as defined in Theorem \ref{TheoremPeriod},  or its $-1$ quadratic twist. Hence, Theorem \ref{TheoremPeriod} provides a formula for the real period of $E$ in this case. If $E(\mathbb{R})$ has two connected components then $E$ is isomorphic to an elliptic curve in Legendre normal form, ${E'}_{\lambda}: y^2=x(x-1)(x-\lambda)$, for some $\lambda \in \mathbb{R} \setminus \lbrace 0, 1 \rbrace$. This is the case covered in \cite{Ro}.

\section{acknowledgements}
The author would like to thank Robert Osburn for his advice during the preparation of this paper and the UCD Ad Astra Research Scholarship program for its financial support.

\end{document}